\documentclass[a4paper, twoside, notitlepage, 12pt]{amsart} 
\usepackage{amsmath,amsfonts,amssymb,amsthm,mathrsfs} 

\usepackage[matrix,arrow]{xy}

\usepackage{verbatim}
\usepackage[latin1]{inputenc}   

\title[The space of generically \'etale families]{The space of generically \'etale families}

\author{David Rydh}
\address{Department of Mathematics \\ KTH \\ Stockholm \\ Sweden}
\email{dary@math.kth.se}

\author{Roy Skjelnes}
\address{Department of Mathematics \\ KTH \\ Stockholm \\ Sweden}
\email{skjelnes@math.kth.se}


\DeclareMathOperator{\id}{id}

\DeclareMathOperator{\ann}{Ker}
\DeclareMathOperator{\Spec}{Spec}
\DeclareMathOperator{\Proj}{Proj}
\DeclareMathOperator{\can}{can}
\DeclareMathOperator{\Tr}{Tr}
\DeclareMathOperator{\determinant}{det}

\newcommand{\ra}{\longrightarrow}

\newcommand{\alt}{\alpha}

\newcommand{\eqbeg}{\begin{equation}}
\newcommand{\eqend}{\end{equation}}
\newcommand{\arbeg}[1]{\begin{array}{#1}}
\newcommand{\arend}{\end{array}}

\newcommand{\p}{\mathfrak{p}}

\newcommand{\n}{\mathfrak{n}}

\newcommand{\Sym}{\mathfrak{S}}

\newcommand{\N}{{P}}
\newcommand{\fpr}{\mathrm{fpr}}

\newcommand{\calA}{\mathscr{A}}
\newcommand{\calR}{\mathscr{R}}
\newcommand{\calZ}{\mathscr{Z}}

\newcommand{\calG}{\mathscr{G}}
\newcommand{\calI}{\mathscr{I}}
\newcommand{\calH}{\mathscr{H}}
\newcommand{\calU}{\mathscr{U}}
\newcommand{\calC}{\mathscr{C}}

\renewcommand{\hom}{\mathrm{Hom}}

\newcommand{\symquot}{\mathrm{Sym}}
\newcommand{\cdiag}{V_X}

\newcommand{\Ten}{\operatorname{T}}
\newcommand{\TS}{\operatorname{TS}}

\newtheorem{thm}[subsection]{Theorem}
\newtheorem{lemma}[subsection]{Lemma}
\newtheorem{cor}[subsection]{Corollary}
\newtheorem{prop}[subsection]{Proposition}

\theoremstyle{definition}

\theoremstyle{remark}
\newtheorem{rem}[subsection]{Remark}

\numberwithin{equation}{subsection}

\newcommand{\calO}{\mathcal{O}}

\begin{document}
\pagenumbering{arabic}

\begin{abstract} 
We construct a space ${\calG}^n_X$ and a rank $n$, generically \'etale family of closed subspaces in a separated ambient space $X$. The constructed pair satisfies a universal property of generically \'etale families of closed subspaces in $X$. This universal property is  derived directly from the construction and does in particular not use the Hilbert scheme. The constructed space ${\calG}^n_X$ is by its universal property canonically identified with a closed subspace of the Hilbert scheme. 
\end{abstract}

\maketitle

\section*{Introduction}

The configuration space of $n$ unordered distinct points on a variety is the complement of the diagonals in the symmetric product of the variety. In the literature there exist several compactifications of the configuration space of $n$ {\em ordered} points, most known is perhaps the Fulton-MacPherson compactification \cite{fultonmacpherson}, but see also \cite{ulyanovcompactification} and the references therein. 

In the relative setting, with a variety replaced with a family of spaces, the space ${\calU}^n_X$ of finite \'etale rank $n$ closed subspaces in $X\ra S$ is analogous to the configuration space of $n$ {\em unordered} points. We are in this article interested in the Hilbert scheme compactification of the space ${\calU}^n_X$. Our work is motivated by the good component construction carried out in \cite{ekedahl&skjelnes}. However, our starting point was to understand the compactification {\it without} using the Hilbert scheme.

When considering this particular compactification we are interested in families of closed subspaces in $X\ra S$ that are flat, finite and almost everywhere \'etale. A flat and finite family $Z\ra S$ is said to be generically \'etale if the open subset of \'etale fibers is schematically dense in the base $S$. Schematically dense in this situation is equivalent with the discriminant being a Cartier divisor.

The property of being generically \'etale is not stable under base change, hence these data do not give a subfunctor of the Hilbert functor of flat families of closed subspaces in $X\ra S$. Still we construct a space ${\calG}^n_X$, for any separated ambient space, and a rank $n$, generically \'etale family ${\calZ}_X\ra {\calG}^n_X$ that satisfies the following universal property: For any $S$-space $T$, and any rank $n$, generically \'etale family $Z \ra T$ of closed subspaces in $X\times_ST$, there exists a unique morphism $f : T\ra {\calG}^n_X$ such that $f^*{\calZ}_X=Z$.

These properties of ${\calG}^n_X$ are derived directly from the construction. In particular they are proved without using the existence of the Hilbert scheme. The fact that the universal properties are easily detected from our construction surprised us, and is the main novelty of the present article.

The relation to the Hilbert functor is clear. The space ${\calU}^n_X$ of \'etale families in $X$ is an open subfunctor of the Hilbert functor, and the schematic closure of ${\calU}^n_X$ is necessarily isomorphic to the constructed space ${\calG}^n_X$. And therefore ${\calG}^n_X$ is also isomorphic to the good component as constructed in  \cite{ekedahl&skjelnes}.

To describe our methods and results in more detail let $X=\Spec(R)$ be an affine scheme over some affine base $S=\Spec (A)$. For any $n$-tuple of elements $x=x_1, \ldots, x_n$ in  $R$ we construct a symmetric tensor $\alt^2(x)$. By restricting the natural map $X^n_S/\Sym_{n-1} \ra X^n_S/\Sym_n$ to the open set where $\alt^2(x)$ does not vanish we get a morphism of affine schemes $\Spec({\calR}) \ra \Spec({\calA})$. The technical focus of the present article is concentrated on the algebra homomorphism ${\calA} \ra {\calR}$. In particular we derive an explicit formula expressing elements in ${\calR}$ as a linear combination of $x_1, \ldots, x_n$.

Using these explicit identities we show that the extension ${\calA}\ra {\calR}$ is \'etale, and that the images of $x_1, \ldots , x_n$ in ${\calR}$ form an ${\calA}$-module basis. Moreover, the pair ${\calA}\ra {\calR}$ is universally determined by these properties (Theorem (\ref{thm:univ.et})). In other words, $\Spec({\calR}) \ra \Spec({\calA})$ parameterizes closed subschemes $Z\subseteq X=\Spec (R)\ra \Spec (A)$ that are \'etale over the base, and such that the images of $x_1, \ldots, x_n$ form a basis for the $A$-module of global sections of $Z$.

A key observation for the representation statement is that the $A$-algebra ${\calA}$ is, in addition to the inverse of $\alt^2(x)$, generated by the traces of endomorphisms of ${\calR}$. The universal property of the extension ${\calA}\ra {\calR}$ is then a matter of specialization.

In the general situation with $X\ra S$ a separated map of algebraic spaces, the universal properties of ${\calZ}_X \ra {\calG}^n_X$ follows from the affine situation combined with certain properties of the blow-up and the existence of the Grothendieck-Deligne norm map.



\subsection*{Structure of the article} In the first section we recall the alternator map and we deduce some explicit formal identities in the tensor algebra. These identities play an important role as they are used to control the relation between two algebras ${\calA}\ra {\calR}$ that we construct. In Section 2 we show that the algebra ${\calR}$ is free as an ${\calA}$-module, and we give a basis. In Section 3 we show that the algebra extension is \'etale, and we prove a universal defining property for the algebra extension. Then in Section 4 we modify our constructed algebras and obtain a new algebra extension ${\calA}_+ \ra {\calR}_+$, and we show its universal defining property.

In Section 5 we consider the global picture, and we let $X\ra S$ be any separated algebraic space. We show how to obtain a space parameterizing closed subspaces of $X$ that are finite, \'etale and of rank $n$ over the base. This gives an independent, abstract, global version of the discussion in Section 3. In Section 6 we construct the space ${\calG}^n_X$ and show that it parameterizes generically \'etale rank $n$ closed subspaces of $X\ra S$, and we show that ${\calG}^n_X$ is the Hilbert scheme compactification of the open subset of \'etale families.

\section{The alternator and symmetric tensors}

In this section we will set up some notation  and establish the identity given in Proposition (\ref{prop:symmetricspan}) that we will use in the next section. 

\subsection{Notation}
Let $R$ be an $A$-algebra. For each non-negative integer $n$ we let $\Ten_A ^nR =R\otimes_A \dots \otimes_A R$ denote the $n$-fold tensor product of $R$ over $A$. The $A$-algebra $\Ten^n_AR$ can be viewed as an $R$-algebra in several different ways. For each $p=1, \ldots , n$, we let $\varphi_p : R \ra \Ten^n_AR$ denote the co-projection on the $p$\textsuperscript{th} factor.

The symmetric group $\Sym _n$ of $n$-letters acts naturally on $\Ten^n_AR$ by permuting the factors. We let $\TS^n_AR:=(\Ten^n_AR)^{\Sym_n}$ denote the $A$-algebra of invariants.

\subsection{The alternator}
\label{subsec:normvector}
For any $n$-tuple of elements $x=x_1, \ldots , x_n$ in $R$ we have the alternating tensor (called norm vector in \cite{ekedahl&skjelnes})
$$ \alt (x_1, \ldots, x_n)=\sum_{\sigma \in \Sym_n}(-1)^{|\sigma|}x_{\sigma(1)}\otimes \dots \otimes x_{\sigma(n)} \quad \text{in} \quad \Ten^n_AR.$$
We will often write $\alt(x)$ instead of $\alt(x_1, \ldots, x_n)$. Let $X$ denote the $(n\times n)$ matrix with coefficients $X_{p,q}:=\varphi_p(x_q)$, where $\varphi_p : R \ra \Ten^n_AR$ are the co-projections. Then we have the determinantal expression of the alternating tensor as 
$$ \alt(x)=\det (X) =\sum_{\sigma \in \Sym_n}(-1)^{|\sigma |}\varphi_{\sigma(1)}(x_1)\cdots \varphi_{\sigma(n)}(x_n).$$
The {\em alternator} is the induced map of  $A$-modules
\eqbeg 
\label{eq:normmap}
\alt : \Ten^n_AR \ra \Ten^n_AR.
\eqend

\begin{prop} 
\label{prop:symlinear}
The alternator \eqref{eq:normmap} is $\TS^n_AR$-linear.
\end{prop}

\begin{proof}
Let $x$ and $y$ be elements of $\Ten^n_AR$, with $y$ a symmetric tensor. We need to check that $\alt({xy})=\alt({x})\cdot y$. As the map $\alt$ is $A$-linear and in particular respects sums, we may assume that $x=x_1\otimes \dots \otimes x_n$. Write $y=\sum_{\gamma}y_{\gamma_1}\otimes \dots \otimes y_{\gamma_n}$. We have
\begin{align*}
\alt({xy}) &= \sum_{\gamma} \alt(x_1y_{\gamma_1}, \ldots, x_ny_{\gamma_n}) \\
&=\sum_{\gamma} \sum_{\sigma \in \Sym_n} (-1)^{|\sigma |}
   x_{\sigma(1)}y_{\gamma_{\sigma(1)}} \otimes \cdots \otimes
   x_{\sigma(n)}y_{\gamma_{\sigma (n)}}                                    \\
&=\sum_{\sigma \in \Sym_n}(-1)^{|\sigma |}
   \big(x_{\sigma(1)}\otimes \cdots \otimes x_{\sigma (n)}\big)\cdot
   \big(\sum_{\gamma}y_{\gamma_{\sigma(1)}}\otimes \cdots \otimes
   y_{\gamma_{\sigma(n)}}\big).
\end{align*}
As $y$ is symmetric we have that $y=\sum_{\gamma}y_{\gamma_{\sigma(1)}}\otimes \cdots \otimes y_{\gamma_{\sigma(n)}}$ for any permutation $\sigma \in \Sym_n$. And consequently that the latter expression obtained above, is ${\alt( x)\cdot y}$.
\end{proof}

\begin{rem} By the usual properties of the determinant we have that the alternator \eqref{eq:normmap} factors as a  map of $A$-modules
$$ \wedge^n_AR \ra \Ten^n_AR.$$
  When two is not a zero divisor in $R$ then $\wedge^n_AR$ is in the natural way a $\TS^n_AR$-module, hence we get an induced $\TS^n_AR$-linear map from the exterior product (see also \cite{laksovthorup_determinantal}). In general, however, the exterior product $\wedge^n_AR$ is not a $\TS^n_AR$-module (see \cite{lundkvist}). We thank D. Laksov for drawing our attention towards these differences.
\end{rem}

\subsection{Alternators of different degrees}
We will compare the alternator maps of different degrees, and in the sequel we
will let $\alt_n$ denote the alternator map whose source is $\Ten^n_AR$.  We let the group $\Sym_{n-1}$ act on $\Ten^n_AR$ by permuting the first $(n-1)$-factors, and let
$$\TS^{n-1,1}_A R=(\Ten^n_A
R)^{\Sym_{n-1}}$$
denote the invariant ring. It is
furthermore convenient to introduce the following notation $\alt_{n-1,1}=\alt_{n-1}\otimes_A \id_R$. An immediate relation between the alternators of different degrees is
\eqbeg \label{identity}
\alpha (x) =\sum_{j=1}^n (-1)^{\mid \tau_{j,n} \mid} \alpha_{n-1,1}(\tau_{j,n}(x)),
\eqend
where $\tau_{j,n} \in \Sym_n$ is the transposition of the factors $j$ and $n$. 

\begin{lemma} The map $\alt_{n-1,1} : \Ten^n_AR \ra \Ten^n_AR$ is $\TS^{n-1,1}_AR$-linear.
\end{lemma}

\begin{proof} By reasoning as in Proposition (\ref{prop:symlinear}) the result follows.
\end{proof}

\begin{prop}
\label{prop:symmetricspan} Let $x=x_1, \ldots, x_n$ be an $n$-tuple of elements in an $A$-algebra $R$. For each integer $i=1, \ldots , n$ we let
$$ x_{[i]} =x_1 \otimes \cdots \otimes x_{i-1} \otimes x_{i+1} \otimes \cdots \otimes x_n \otimes 1 \quad \text{in}\quad \Ten^n_AR$$
denote the $n$-tensor we get by removing the $i$'th factor and inserting $1$ on the last factor. We  have, for any $y\in \TS^{n-1,1}_AR$, the identity
$$ \alt(x)y= \sum_{i=1}^n(-1)^{n-i}\alt\big( x_{[i]} y\big)\varphi_n(x_i)\quad \text{in} \quad \Ten^n_AR.$$
\end{prop}

\begin{proof} Using the identity (\ref{identity}) we get that $\alpha(x_{[i]}y)$ can be written as $\sum_{j=1}^n(-1)^{\mid \tau_{j,n}\mid}\alpha_{n-1,1}(\tau_{j,n}(x_{[i]}y))$. The element $\varphi_n(x_i)=1\otimes \cdots \otimes 1 \otimes x_i$ is $\Sym_{n-1}$-invariant. And as $\alpha_{n-1,1}$ is $\TS^{n-1,1}_AR$-linear we have that the sum appearing in the proposition is the sum of
\eqbeg
\label{eq:expr1}
\sum_{i=1}^n(-1)^{n-i}\alt_{n-1,1}\big( x_{[i]}\varphi_n(x_i) y\big)
\eqend
and
\eqbeg
\label{eq:expr2}
\sum_{i=1}^n\sum_{j=1}^{n-1} (-1)^{n-i+1}
\alt_{n-1,1}\big(\tau_{j,n}(x_{[i]}y)\varphi_n(x_i)\big).
\eqend
In the first sum (\ref{eq:expr1}) we can take out the $\Sym_{n-1}$-invariant element $y$ by the $\TS^{n-1,1}_AR$-linearity of $\alpha_{n-1,1}$. It is readily checked that 
$$\sum_{i=1}^n(-1)^{n-i}\alpha_{n-1,1}(x_{[i]}\varphi_n(x_i))=\alpha (x).$$ Thus (\ref{eq:expr1}) equals $\alpha(x)y$ and we need to check that the remaining expression (\ref{eq:expr2}) is zero. The summation in (\ref{eq:expr2}) runs over the set of pairs  $(i,j)$ with $i=1, \ldots, n$ and $j=1, \ldots , n-1$, which is the disjoint sum of ${\calC}_{>}=\{(i,j) \mid i>j\}$ and ${\calC}_{\leq}=\{(i,j) \mid i\leq j\}$. For fixed $(i,j) \in {\calC}_{>}$ we have the corresponding summand 
\eqbeg \label{sum>}
(-1)^{n-i+1}\alpha_{n-1,1}(\tau_{j,n}(x_{[i]}y)\varphi_n(x_i))
\eqend
in $(\ref{eq:expr2})$. We claim that this summand is canceled by the corresponding summand over $(j,i-1) \in {\calC}_{\leq}$. Let $\sigma^q_p \in \Sym_n$, for any integers $p\leq q\leq n$, denote the cyclic, increasing permutation of the factors $p, \ldots, q$. Since $j<i$ we have
\eqbeg \label{id2}
 \tau_{j,n}(x_{[i]}y)\varphi_n(x_i) =\sigma^{i-1}_j\big(\tau_{i-1,n}(x_{[j]})\varphi_n(x_j)\big)\tau_{j,n}(y).
\eqend
Moreover,  $\sigma^{i-1}_j \circ \tau_{i-1,n}=\tau_{j,n}\circ \sigma^{i-1}_j$, and as $\sigma^{i-1}_j(y)=y$ we can write (\ref{id2}) as
$\sigma^{i-1}_j \big(\tau_{i-1,n}(x_{[j]}y)\varphi_n(x_j)\big)$.
It then follows that (\ref{sum>}) equals 
$$(-1)^{n-i+1+\mid \sigma^{i-1}_j \mid } \alpha_{n-1,1}(\tau_{i-1,n}(x_{[j]}y)\varphi_n(x_j)).$$
As $\mid \sigma^{i-1}_j \mid =i-1-j$ we have that the summand (\ref{sum>}) is canceled by the corresponding summand over $(j,i-1) \in {\calC}_{\leq}$.  Consequently (\ref{eq:expr2}) is zero and we have proven the proposition.
\end{proof}

\subsection{Linear span}
\label{subsec:linspan} Let $z=\sum_{i=1}^na_ix_i$ in $R$, with scalars $a_1, \ldots , a_n$ in $A$. Then, as the alternator $\alt : \Ten^n_AR \ra \Ten^n_AR$ is multi-linear and alternating, we have
\eqbeg
\label{eq:coefficient} \alt(x_1, \ldots, x_{i-1}, z, x_{i+1}, \ldots, x_n) = a_i \cdot \alt(x) \quad \text{in}\quad \Ten^n_AR.
\eqend
In particular these identities hold for any $z\in R$ if the elements $x_1, \ldots , x_n$ form an $A$-module basis of $R$. However these identities always hold formally in $\Ten^n_AR$. Indeed, as a special case of our previous result we have

\begin{cor}
\label{cor:Rspan}
Let $x=x_1, \ldots, x_n$ be an $n$-tuple of elements in $R$. For any $z\in R$ we have the identity
$$\alt(x)\varphi_n(z)=\sum_{i=1}^n \alt(x_1, \ldots, x_{i-1}, z, x_{i+1}, \ldots , x_n)\varphi_n(x_i)$$
in the tensor algebra $\Ten^n_AR$.
\end{cor}

\begin{proof} Specializing the proposition with $y=\varphi_n(z) \in \TS^{n-1,1}_AR$ gives $\alt(x)\varphi_n(z) = \sum_{i=1}^n(-1)^{n-i}\alt(x_{[i]}\varphi_n(z))\cdot \varphi_n(x_i)$. The result then follows since $(-1)^{n-i}\alt(x_{[i]}\varphi_n(z))=\alt(x_1, \ldots, x_{i-1},z,x_{i+1}, \ldots, x_n)$.
\end{proof}

\begin{rem} There is a canonical map $\TS^{n-1}_A R\otimes_A R\to \TS^{n-1,1}_A R$ but this
is not always an isomorphism in positive characteristic.  In particular the inclusion
$\TS^n_A R\ra \Ten^n_A R$ that factors through the inclusion
$\TS^{n-1,1}_A R\ra \Ten^n_A R$, does not always factorize through
$\TS^{n-1}_A R\otimes_A R$. We thank T. Ekedahl for pointing out this and thereby correcting a mistake in an earlier version of this article.
\end{rem}

\section{Linear solution spaces}

Let $x_1, \ldots, x_n$ be a fixed $n$-tuple of elements in $R$. We will show how to obtain an $R$-algebra $\calR$ and an $A$-subalgebra $\calA\subseteq \calR$ such that the elements $x_1, \ldots ,x_n$ in the $R$-algebra ${\calR}$ form an $\calA$-module basis. The universality of the constructed pair will be established in the section following this one.

\subsection{Localization of the square of an alternating tensor}
\label{subsec:norm^2}
\label{subsec:baseringA} Let $x=x_1, \ldots, x_n$ be an $n$-tuple of elements in $R$, and let $\alt(x)$ denote the tensor described in \ref{subsec:normvector}. It is clear that a permutation $\sigma \in \Sym_n$ sends the tensor $\alt(x) \in \Ten^n_AR$ to $(-1)^{|\sigma|}\alt(x)$. Consequently, if $y=y_1, \ldots, y_n$ is another $n$-tuple of elements in $R$ then $\alt(x)\alt(y)$ is a symmetric tensor, that is $\alt(x)\alt(y) \in \TS^n_AR$. In particular $\alt^2(x)$, the square of the alternating tensor, is symmetric. 

From Corollary (\ref{cor:Rspan}) we see that the element $\alt(x)$ has to be inverted in order to solve the linear equations $z=\sum_{i=1}^na_ix_i$ in $\Ten^n_AR$. We therefore introduce the notation 
$$\calA = \calA(\alt^2(x))= \TS^n_AR[\alt^2(x)^{-1}]$$
for the localization of $\TS^n_AR$ in the symmetric tensor $\alt^2(x)$.

\begin{lemma}
\label{lemma:linindp} The elements $\varphi_n(x_1), \ldots, \varphi_n(x_n)$ in the localized ring $\Ten^n_AR[\alt(x)^{-1}]$ are linearly independent over ${\calA}$.
\end{lemma}

\begin{proof} An element $\sum_{i=1}^na_i'\varphi_n(x_i)$ in $\Ten^n_AR[\alt(x)^{-1}]$ with $a_i' \in {\calA}$ can be written as $\alt^{-2p}(x)\sum_{i=1}^na_i\varphi_n(x_i)$, with $a_i \in \TS^n_AR$, for some integer $p$. Assume therefore that we have a relation
\eqbeg
\label{eq:lindep}
 a_1\varphi_n(x_1)+\cdots + a_n\varphi_n(x_n)=0 
\eqend
in $\Ten^n_AR$, with symmetric tensors $a_1, \ldots , a_n$ in $\TS^n_AR$. A 
 permutation $\sigma \in \Sym_n$ will transform the equation \eqref{eq:lindep}, to the identity 
$$ a_1\varphi_{\sigma(n)}(x_1) +\cdots + a_n\varphi_{\sigma(n)}(x_n)=0.$$
In particular we can replace $\varphi_n$ in \eqref{eq:lindep} with $\varphi_j$, for any $j=1, \ldots, n$. We then obtain the matrix equation
$$ X \cdot A =0,$$
where $X$ is the $(n\times n)$-matrix with coefficients $X_{p,q}=\varphi_p(x_q)$, and $A$ is the $(n\times 1)$-matrix with coefficients $a_1, \ldots, a_n$. The coefficients of the matrices are in $\Ten^n_AR$. The determinant of the matrix $X$ is $\alt(x_1, \ldots , x_n)$, and it follows that the coefficients $a_1, \ldots , a_n$ are all zero in $\Ten^n_AR[\alt(x)^{-1}]$.
\end{proof}

\subsection{The linear solution spaces}
\label{subsec:linspaces} We have the inclusion of rings
\eqbeg
\label{eq:modulemap}
\xymatrix{
\TS^n_AR \ar[r]^-{i} & \TS^{n-1,1}_AR \ar[r]^-{j} & \Ten^n_AR.
}
\eqend
We have the $n$-tuple $x=x_1, \ldots , x_n$ fixed, and we localize the above sequence of rings with respect to the symmetric tensor $\alt^2(x) \in \TS^n_AR$. We have earlier (\ref{subsec:baseringA}) introduced the notation
$\calA:= \TS^n_AR[\alt^2(x)^{-1}]$, and now we introduce 
\eqbeg
\label{eq:linsolspace}
\calR = \calR(\alt^2(x))=\TS^{n-1,1}_AR [\alt^2(x)^{-1}]
\eqend
for the localization of $\TS^{n-1,1}_AR$ with respect to the element $\alt^2(x) \in \TS^n_AR$. The $R$-algebra structure on $\calR$ is the one induced from the co-projection map $\varphi_n : R \ra \Ten^n_AR$ on the last factor. We will also denote the structure map with $\varphi_n : R \ra {\calR}$. Moreover, as the first map $i$ of \eqref{eq:modulemap} is injective, we have that $\calA \subseteq \calR$ is a ring extension.
 
\begin{thm}
\label{thm:basis}
Let $x=x_1, \ldots, x_n$ be an $n$-tuple of elements in an $A$-algebra $R$, and define the ring extension $\calA \subseteq \calR$ as in (\ref{subsec:linspaces}). Then we have that the elements $\varphi_n(x_1), \ldots ,\varphi_n(x_n)$ in the $R$-algebra $\calR$ form an ${\calA}$-module basis.
\end{thm}

\begin{proof} We obtain a map ${\calR}\ra \Ten^n_AR[\alt(x)^{-1}]$ by localizing the canonical map $j : \TS^{n-1,1}_AR \ra \Ten^n_AR$ of \eqref{eq:modulemap}. Linear independence follows from Lemma (\ref{lemma:linindp}). It remains to see that every tensor $y\in \TS^{n-1,1}_A R$ is in the linear
span of the elements $\varphi_n(x_1), \ldots, \varphi_n(x_n)$. From Proposition
(\ref{prop:symmetricspan}) we obtain the identity
$$ \alt^2(x)y =\sum_{i=1}^n (-1)^{n-i}\alt(x) \alt( x_{[i]} y) \varphi_n(x_i).$$
The product of two alternating tensors is a symmetric tensor, hence when inverting $\alt^2(x)$ then $y$ is in the linear span of $\varphi_n(x_1), \ldots, \varphi_n(x_n)$.
\end{proof}


\section{Universal properties of ${\calA}\ra {\calR}$}

In this section we will show that the extension ${\calA}\subseteq {\calR}$ constructed in the previous section, has a universal property. We will show that the pair represent the functor of pairs of \'etale extensions $B\ra E$, where $E$ is an $R$-algebra and where the images of $x_1, \ldots , x_n$ in $E$ form a $B$-module basis.

\subsection{Trace map}
Recall that the extension ${\calA} \subseteq {\calR}$ we have constructed fits into the commutative diagram of $A$-algebras and $A$-algebra homomorphisms
%
$$
\xymatrix{
R \ar[r] & R\otimes_A{\calA} \ar[r] & {\calR}\\
A \ar[u] \ar[r] & {\calA} \ar[u] & 
}
$$
The composition of the two upper horizontal morphisms in the diagram above is $\varphi_n : R \ra {\calR}$. As ${\calR}$ is free of rank $n$ as an ${\calA}$-module, we have the ${\calA}$-linear trace map ${\calR} \ra {\calA}$.

\subsection{Discriminant} Let $B\ra E$ be a homomorphism of rings, where $E$ is free of rank $n$ as a $B$-module. To a given $B$-module basis $e_1, \ldots, e_n$ of $E$ we define the discriminant $d_{E}\in B$ as the determinant of the $(n\times n)$-matrix whose coefficient $c_{i,j}$ is the trace of the $B$-linear map $z\mapsto z \cdot e_ie_j$. It is easily seen that the discriminant $d_{E}$ depends on the choice of basis of $E$, however the ideal $D_{E/B}\subseteq B$ generated by the discriminant does not. If $B\ra E$ is \emph{locally free} it is thus clear that there is a locally principal ideal $D_{E/B}\subseteq B$, locally given by the discriminant.

\begin{lemma} 
\label{lemma:Trace}
For any element $z\in R$ we have that the trace of the ${\calA}$-linear endomorphism $e\mapsto \varphi_n(z)e$ on ${\calR}$ is $\varphi_1(z)+\cdots +\varphi_n(z)$.
\end{lemma}

\begin{proof} The elements $\varphi_n(x_1), \ldots , \varphi_n(x_n)$ form, by Theorem (\ref{thm:basis}) an ${\calA}$-module basis of ${\calR}$. The action of $\varphi_n(z)$ on the fixed basis element $\varphi_n(x_k)$ is $\varphi_n(zx_k)$. Using Corollary (\ref{cor:Rspan}) we have that $\varphi_n(zx_k)$ is expressed as the sum
$$ \sum_{i=1}^n \frac{\alt (x_1, \ldots, x_{i-1}, zx_k, x_{i+1}, \ldots , x_n)\alt(x)}{\alt^2(x)} \varphi_n(x_i).$$ 
In particular we see that $\alt(x_1, \ldots, x_{k-1}, zx_k, x_{k+1}, \ldots, x_n)\alt(x)^{-1}$ is the $k$\textsuperscript{th} component of $\varphi_n(zx_k)$, and consequently that the trace of the endomorphism $e\mapsto \varphi_n(z)e$ is
\eqbeg
\label{eq:tracesum}
\sum_{k=1}^n \frac{\alt(x_1, \ldots, x_{k-1}, zx_k, x_{k+1}, \ldots, x_n)}{\alt(x)}.
\eqend
Now using the fact that the map $\alt$ is $\TS^n_AR$-linear (\ref{prop:symlinear}) we get that 
$$ \alt(x)\sum_{k=1}^n\varphi_k(z)=\sum_{k=1}^n {\alt(x_1, \ldots, x_{k-1}, zx_k, x_{k+1}, \ldots, x_n)}. $$
Applying that identity to the sum \eqref{eq:tracesum} proves the lemma.
\end{proof}

\subsection{Polarized power sums}
For any element $z\in R$ the symmetric tensor
$$ \p (z) := \varphi_1(z)+\cdots +\varphi_n(z) \in \TS^n_AR,$$
is often referred to as polarized power sum \cite{weyl_invariants}. One can check that the symmetric tensor $\alt^2(x)$ is a sum of products of polarized power sums: In fact if $x=x_1, \ldots , x_n$ and $y=y_1, \ldots , y_n$ are two $n$-tuples of elements in $R$, we have~\cite[Lemma 2.3]{ekedahl&skjelnes}
\eqbeg
\label{eq:traceexp}
\alt(x)\alt(y) = \det (\p(x_iy_j)).
\eqend

\begin{prop}
\label{prop:discriminant}
 The discriminant of the extension ${\calA} \ra {\calR}$ is $\alt^2(x)$. In particular we have that ${\calA} \ra {\calR}$ is \'etale.
\end{prop}
\begin{proof} We have by \eqref{eq:traceexp} that $\alt^2(x)$ is the determinant of the $(n\times n)$-matrix whose coefficient $(i,j)$ equals $\p(x_ix_j)$. By definition we have that $\p(x_ix_j)=\varphi_1(x_ix_j)+\cdots +\varphi_n(x_ix_j)$, which by Lemma (\ref{lemma:Trace}) equals the trace of the endomorphism $e\mapsto \varphi_n(x_ix_j)\cdot e$. As $\varphi_n(x_1), \ldots ,\varphi_n(x_n)$ is an ${\calA}$-module basis of ${\calR}$ (Theorem \ref{thm:basis}), we have that $\alt^2(x)$ is the discriminant of the extension ${\calA}\ra {\calR}$. It then follows ~\cite[18.2]{egaIV} that the extension is \'etale.
\end{proof}

\subsection{Algebra of polarized powers}
We let $\N \subseteq \TS^n_AR$ be the $A$-subalgebra generated by the polarized power sums; 
$$ \N =A[ \p(z) ]_{z \in R }.$$
By \eqref{eq:traceexp} we have that $\alt^2(x) \in \N$.

\begin{prop}
\label{prop:polpowalgebra} The inclusion of $A$-algebras $\N \subseteq \TS^n_AR$ becomes an isomorphism after localization with respect to $\alt^2(x)$.
\end{prop}

\begin{proof} Let $x=x_1\otimes \cdots \otimes x_n$. If $y\in \TS^n_AR$ is a symmetric tensor we have by Proposition (\ref{prop:symlinear})
$$ \alt^2(x) \cdot y = \alt(x) \cdot \alt( xy),$$
and therefore $y=\alt(x)\alt( xy)/\alt^2(x)$ in $\TS^n_AR[\alt^2(x)^{-1}]={\calA}$. The symmetric tensor $y$ is a sum of tensors $\sum_{\gamma}y_{\gamma}$, where each summand $y_{\gamma}$ is  on the form $y_{\gamma_1}\otimes  \cdots \otimes y_{\gamma_n}$. Consequently $\alt(xy)$ is the sum $\sum_{\gamma}\alt(xy_{\gamma})$. From \eqref{eq:traceexp} it follows that $\alt(x)\alt(xy)$ is a sum of products of polarized powers. In other words, the symmetric tensor $y\in \N[\alt^2(x)^{-1}]$.
\end{proof}

\begin{rem} When the base ring $A$ contains the field of rationals then we have that the $A$-algebra of symmetric tensors $\TS^n_AR$ is generated by its polarized power sums (see \cite{weyl_invariants}). It is moreover known that the polarized power sums do not always generate $\TS^n_AR$ in positive characteristic. We will later see (\ref{rem:diagonals}) that the support of $\alt^2(x)$ in $\Spec (\TS^n_AR)$, when running through all $n$-tuples $x$, is precisely the diagonals. What our proposition above therefore says is that the polarized power sums always generate the ring of invariants - as long as we stay away from the diagonals.
\end{rem}

\begin{prop}
\label{prop:whenRisfree}
Assume that the elements $x_1, \ldots , x_n$ in $R$ form an $A$-module basis. Then the discriminant $d_R \in A$ is mapped to $\alt^2(x)$ by the structure map $A\ra \calA$. The induced map $A_{d_R} \ra {\calA}$ is an isomorphism, and the $A$-algebra homomorphism $\varphi_n : R\otimes_AA_{d_R} \ra {\calR}$ is an isomorphism.
\end{prop}

\begin{proof} We note that the canonical $A$-algebra morphism ${\calA}\otimes_AR \ra {\calR}$ is an isomorphism. Indeed, bijectivity follows as both ${\calA}$-modules have the same basis: The elements $1\otimes x_1, \ldots , 1\otimes x_n$ is a basis of ${\calA}\otimes_AR$ by assumption, and these elements are mapped to $\varphi_n(x_1), \ldots, \varphi_n(x_n)$ which form a basis of ${\calR}$ by Theorem (\ref{thm:basis}). In particular we have, for any $z\in R$, the following identity of traces:
\eqbeg
\label{eq:trace&basechange}
\Tr_A (e\mapsto z e) \otimes 1 = \Tr_{\calA}(e\mapsto \varphi_n(z) e).
\eqend
By Lemma (\ref{lemma:Trace}) we then have that $\p(z)=\Tr_A(e\mapsto ze)\otimes 1$, and by Proposition (\ref{prop:discriminant}) we have that the discriminant $d_R\in A$ of the extension $A\ra R$ is mapped to $\alt^2(x) \in {\calA}$ by the structure map $A\ra {\calA}$. It now follows from \eqref{eq:trace&basechange} that the induced map $A_{d_R} \ra {\calA}$ is surjective. Indeed, by Proposition (\ref{prop:polpowalgebra}) the $A$-algebra ${\calA}$ is generated by the polarized powers $\p(z)$ and the inverse of $\alt^2(x)$.

It remains to the see injectivity. As our morphism $A_{d_R} \ra {\calA}$ is the localization of the map $A\ra \TS^n_AR$, it suffices to see that the map $A \ra \Ten^n_AR$ is injective. Since the algebra $R$ is a free $A$-module, it follows that $\Ten^n_AR$ is free, hence $A\ra \Ten^n_AR$ is injective.
\end{proof}

\subsection{A commutative diagram}
\label{subsec:commdiag}
Let $B\ra E$ be a homomorphism of $A$-algebras. If $f : R\ra E$ is an $A$-algebra homomorphism, we obtain a natural $A$-algebra homomorphism of tensor products $\Ten^n_AR \ra \Ten^n_BE$, taking $z_1\otimes \cdots \otimes z_n$ to $f(z_1)\otimes \cdots \otimes f(z_n)$. Consequently we get induced $A$-algebra homomorphisms $f_n : \TS^n_AR \ra \TS^n_BE$ and $f_{n-1,1} : \TS^{n-1,1}_AR \ra \TS^{n-1,1}_BE$ such that the following diagram of $A$-algebras
\eqbeg
\label{eq:commdiag}
\vcenter{\xymatrix{
\TS^n_AR \ar[r]^{f_n} \ar[d] & \TS^n_BE \ar[d] \\
\TS^{n-1,1}_AR \ar[r]^{f_{n-1,1}}& \TS^{n-1,1}_BE
}}
\eqend
commutes.

\begin{thm}
\label{thm:univ.et} Let $B\ra E$ be an \'etale extension of $A$-algebras, and let $f : R \ra E$ be an $A$-algebra homomorphism. Assume that the elements $f(x_1), \ldots , f(x_n)$ form a $B$-module basis of $E$. Then there is a unique $A$-algebra homomorphism
$$ \n_{E/B} : {\calA} \ra B,$$
such that ${\calR}\otimes_{\calA} B=E$ as quotients of $R\otimes_AB$. The homomorphism $\n_{E/B}$ is determined by sending the polarized power sum $\p(z)$ to $\Tr_B(e\mapsto f(z)e)$, the trace of the multiplication map by $f(z)$, for every $z\in R$.
\end{thm}

\begin{proof}  We have the symmetric tensor  $\alt^2(f(x_1), \ldots, f(x_n))=\alt^2(f(x))$ in $\TS^n_BE$.  Let ${\calA}_E$ and ${\calR}_E$ denote the localization with respect to $\alt^2(f(x))$ in  $\TS^n_BE$ and of $\TS^{n-1,1}_BE$, respectively. The natural morphism $f_n : \TS^n_AR \ra \TS^n_BE$ takes $\alt^2(x)$ to $\alt^2(f(x))$. Then by localizing  the commutative diagram \eqref{eq:commdiag} with respect to $\alt^2(x)$ we obtain an ${\calA}$-algebra homomorphism
\eqbeg
\label{eq:inducedmap}
{\calR}\otimes_{\calA}{\calA}_E \ra {\calR}_E.
\eqend

The morphism \eqref{eq:inducedmap} takes the basis $\varphi_n(x_1)\otimes 1, \ldots , \varphi_n(x_n) \otimes 1$ of ${\calR}\otimes_{\calA}{\calA}_E$ to the the elements $\varphi_n(f(x_1)), \ldots , \varphi_n(f(x_n))$. By Theorem (\ref{thm:basis}) we have that $\varphi_n(f(x_1)), \ldots , \varphi_n(f(x_n))$ form an ${\calA}_E$-module basis of ${\calR}_E$. Consequently the ${\calA}$-algebra homomorphism \eqref{eq:inducedmap} is a surjective map of free rank $n$ ${\calA}_E$-modules, hence an isomorphism.

By Proposition (\ref{prop:whenRisfree}) we have a natural identification $B={\calA}_E$, and $E={\calR}_E$, and consequently we have obtained an $A$-algebra homomorphism $\n_{E/B} : {\calA} \ra B$ such that ${\calR}\otimes_{\calA}B=E$. 

It remains to see that the morphism ${\n_{E/B}}$ is unique, and that it maps polarized powers $\p(z)$ to the trace of the multiplication map on $E$. Note that since ${\calR}\otimes_{\calA}B=E$ we have the identity of traces 
$ \Tr_{\calA}(e\mapsto \varphi_n(z)e)\otimes 1 =\Tr_B(e\mapsto f(z)e)$,
for any element $z\in R$. Combining that with Lemma (\ref{lemma:Trace}) yield the identities 
$$ \n_{E/B}(\p(z))=\Tr_{\calA}(e\mapsto \varphi_n(z)e)\otimes 1 = \Tr_B(e\mapsto f(z)e).$$
From where we conclude that the morphism $\n_{E/B} : {\calA} \ra B$ has the announced action on $\p(z)$, and since $\p(z)$ generates ${\calA}$ (Proposition \ref{prop:polpowalgebra}) the morphism is unique.
\end{proof}

\subsection{The universal \'etale family, local description}
For each $n$-tuple of elements $x=x_1, \ldots, x_n$ in $R$ we have the open immersion of affine schemes $\Spec (\calA (\alt^2(x)) \subseteq \Spec (\TS^n_AR)$. Let $\bf{x}$ denote the set of all $n$-tuples of elements in $R$, and consider the union
$$ \calU_{R}=\bigcup_{x\in{\bf x}}\Spec (\calA (\alt^2(x)) \subseteq \Spec (\TS^n_AR).$$
By construction we have a map $\calZ_{R}=\cup_{x\in{\bf x}}\Spec (\calR (\alt^2(x))\ra \calU_{R}$.

\begin{cor} 
\label{cor:etalaffine} The family $\calZ_{R}\ra \calU_{R}$ is \'etale of rank $n$, and the natural map $\calZ_{R}\ra \Spec (R)\times_{\Spec(A)}\calU_{R}$ is a closed immersion. Furthermore, let $B$ be an $A$-algebra, and let $E$ be an algebra quotient of $R\otimes_AB$ which is \'etale of rank $n$ as a $B$-algebra. Then there exists a unique morphism of schemes $\n : \Spec (B) \ra \calU_{R}$ such that the pull-back $\n^{*}\calZ_{R}=\Spec (E)$ as closed subschemes of $\Spec (R\otimes_AB)$.
\end{cor}

\begin{proof} It is clear from the discussion in this section that the family $\calZ_{R}\ra \calU_{R}$ is a closed subscheme of $\Spec (R)\times_{\Spec (A)}\calU_{R}$ and \'etale of rank $n$ over $\calU_{R}$. We will show its universal properties, and we begin with uniqueness.

Assume that there are two morphisms
$\n_1,\n_2 : \Spec(B) \ra \calU_R$ such that the
two pull-backs $\n_1^{*}\calZ_R=\n_2^{*}\calZ_R$ coincide with $\Spec (E)$. For any $n$-tuple $x$ of elements in $R$,  write $\calU_{R,x}=\Spec\bigl(\calA(\alt^2(x))\bigr)$. We define the open
subsets $T_{x,y}=\n_1^{-1}(\calU_{R,x}\bigr)\cap\n_2^{-1}(\calU_{R,y})$, for any two $n$-tuples $x,y \in \bf{x}$ of elements in $R$. Then
$\{T_{x,y}\}_{x,y\in \bf{x}}$ is an open cover of $\Spec(B)$. Moreover $T_{x,y}$
is affine as $\calU_{R,x}$ is affine. We let
$T_{x,y}=\Spec(B_{x,y})$. The images in $E_{x,y}=E\otimes_B B_{x,y}$ of either
the elements $x=x_1, \ldots, x_n$ or the elements $y=y_1, \ldots, y_n$ form a $B_{x,y}$-module basis.

Let $\n_x : \calA(\alt^2(x)) \ra B_{x,y}$ and
$\n_y : \calA(\alt^2(y)) \ra B_{x,y}$ be the homomorphisms corresponding to
the morphisms of schemes $\n_1|_{T_{x,y}}$ and $\n_2|_{T_{x,y}}$.
By Theorem (\ref{thm:univ.et}) $\n_x$ and $\n_y$ sends polarized power sums
onto traces. Thus, by \eqref{eq:traceexp}, both $\n_x$ and $\n_y$ sends
$\alt^2(x)$ and $\alt^2(y)$ onto a generator of the discriminant
ideal $D_{E_{x,y}}\subseteq B_{x,y}$. By the \'etaleness assumption of $E$
we thus have that both $\n_x$ and $\n_y$ factors through
$\TS^n_AR[\alt^2(x)^{-1}\alt^2(y)^{-1}]$. By the uniqueness of $\n_x$ and $\n_y$
it follows that these factorizations are equal. Thus $\n_1=\n_2$.

Now, let $B$ and $E$ be as in the corollary. Then there is an affine open
covering $\{\Spec(B_{\gamma})\}_{\gamma}$ of $\Spec (B)$ such that for each
index $\gamma$ there are $n$ elements $x_{\gamma}=x_{\gamma_1}, \ldots , x_{\gamma_n}$ in $R$ whose images in $E_{\gamma}=E\otimes_BB_{\gamma}$ form a $B_{\gamma}$-module basis. By Theorem (\ref{thm:univ.et}) we get, for each
$\gamma$, a morphism
$$ \n_{\gamma} : \Spec(B_{\gamma}) \ra \calU_{R,x_\gamma} \subseteq \calU_{R}$$
such that the family $\calZ_{R}$ is pulled back to
$\Spec (E_{\gamma})$. The restrictions of $\n_{\gamma}$ and $\n_{\gamma'}$ to
the intersection $\Spec(B_\gamma)\cap\Spec(B_{\gamma'})$ coincide by the
uniqueness shown above. Thus the maps $\{\n_{\gamma}\}$ glue to a morphism
$\n : \Spec(B) \ra \calU_{R}$ with the required property.
\end{proof}

\section{Generically \'etale families}

By some minor modifications of the results in the previous section we will obtain a representing pair ${\calA}_+ \ra {\calR}_+$ for the functor of generically \'etale families. We use the notation introduced in the preceding sections, and in particular we have the $n$-tuple $x=x_1, \ldots, x_n $ of elements in $R$ fixed.

\subsection{A canonical ideal}
\label{subsec:normideal}
We have (Proposition \ref{prop:symlinear}) that the alternator map $\alt : \Ten^n_AR \ra \Ten^n_AR$ is $\TS^n_AR$-linear, and since $\Ten^n_AR$ has a product structure we obtain a $\TS^n_AR$-module map
$$ \alt \times \alt : \Ten^n_AR \otimes_{\TS^n_AR}\Ten^n_AR \ra \Ten^n_AR,$$
taking two elements $y$ and $z$ of $\Ten^n_AR$ to $\alt(y)\alt(z)$. From (\ref{subsec:norm^2}) we have that $\alt(y)\alt(z)$ is invariant under the $\Sym_n$-action. Consequently the $\TS^n_AR$-module given as the image of $\alt\times \alt $ is an ideal of $\TS^n_AR$, and this ideal we will denote by $I_{R}$. We will refer to it as the {\em canonical ideal}.

\begin{rem} There is a natural homomorphism of rings from the divided powers ring $\Gamma^n_AR$ to the invariant ring $\TS^n_AR$. Under this map the canonical ideal $I_R$ is the image of the norm ideal $I\subseteq \Gamma^nA_R$ defined in (\cite{ekedahl&skjelnes}).
\end{rem}

\subsection{The subscheme defined by the canonical ideal}
\label{subsec:nu^2(x)cover} It is clear from the definition of the canonical ideal $I_{R}\subseteq \TS^n_AR$  that it is generated by elements of the form $\alt(y)\alt(z)$, with $n$-tuples $y=y_1, \ldots, y_n$ and $z=z_1, \ldots , z_n$ of elements in $R$. Furthermore, if we let $D(\alt^2(x))\subseteq \Spec (\TS^n_AR)$ denote the open subset where $\alt^2(x)$ does not vanish, then we have
$$ D(\alt(y)\alt(z))=D(\alt^2(y))\cap D(\alt^2(z)) \quad \text{ in }\quad \Spec (\TS^n_AR).$$
Let $\Delta \subseteq \Spec (\TS^n_AR)$ denote the closed subscheme defined by $I_{R}$. We then have, using the notation of the preceding section, that
\eqbeg
\label{eq:deltacomplement}
\calU_{R} =\bigcup_{x\in {\bf x}} \Spec (\calA({\alt^2(x)}))=\Spec (\TS^n_AR) \setminus \Delta,
\eqend
where ${\bf x}$ is the set of all $n$-tuples of elements in $R$. Let 
$$\psi : \Spec (\TS^{n-1,1}_AR) \ra \Spec (\TS^n_AR)$$
denote the morphism corresponding to the canonical ring homomorphism $i : \TS^n_AR \ra \TS^{n-1,1}_AR$ \eqref{eq:modulemap}. We then have, with the notation of Corollary (\ref{cor:etalaffine}), that the family ${\calZ}_{R}=\psi^{-1}({\calU}_{R})$.

\subsection{The blow-up algebras}
The blow-up algebra $\oplus_{m\geq 0}I^m_{R}$ of the canonical ideal $I_{R}\subseteq \TS^n_AR$ is a graded ring, and we let
$$ {\calA}_+ = {\calA}_+(\alt^2(x))=\big(\oplus_{m\geq 0}I_{R}^m\big)_{(\alt^2(x))}$$
denote the degree zero part of the localization at the degree one element $\alt^2(x) \in I_{R}$ of the blow-up algebra $\oplus_{m\geq 0}I^m_{R}$. The natural inclusion of rings $ \TS^n_AR \ra \TS^{n-1,1}_AR$ \eqref{eq:modulemap} induces a graded ring homomorphism between their blow-up algebras, and then an induced ring homomorphism
\eqbeg
\xymatrix{
{\calA}_+\ar[r] & {\calR}_+:=\big(\oplus_{m\geq 0}I^m_{R}\TS^{n-1,1}_AR\big)_{(\alt^2(x))}.
}
\eqend
The $A$-algebra homomorphism $\varphi_n : R\ra \TS^{n-1,1}_AR$ gives an $R$-algebra structure on ${\calR}_+$.

\begin{prop}
\label{prop:specialcase} 
The ${\calA}_+$-algebra ${\calR}_+$ is free as an ${\calA}_+$-module, and the elements $\varphi_n(x_1), \ldots , \varphi_n(x_n)$ form a basis. In particular we have that the $A$-algebra homomorphism $R\otimes_A{\calA}_+ \ra {\calR}_+$ is surjective. Furthermore, the discriminant of the extension ${\calA}_+ \ra {\calR}_+$ is a non-zero divisor in ${\calA}_+$.
\end{prop}

\begin{proof} Linear independence of $\varphi_n(x_1), \ldots , \varphi_n(x_n)$ follows from Theorem (\ref{thm:basis}) as the localization of ${\calA}_+$ in the element $\alt^2(x) \in \TS^n_AR$ is injective, and we have $({\calA}_+)_{\alt^2(x)}={\calA}$. From Proposition (\ref{prop:symmetricspan}) elements of the form $y\in\TS^{n-1,1}_AR$ are in the linear span of $\varphi_n(x_1), \ldots , \varphi_n(x_n)$, and consequently they form a basis.

It remains to see that the discriminant of the extension ${\calA}_+\ra {\calR}_+$ is a non-zero divisor. As the localization map ${\calA}_+\ra {\calA}$ is injective, and ${\calR}_+\otimes_{{\calA}_+}{\calA}={\calR}$ we need only to see that the discriminant of ${\calA}\ra {\calR}$ is a non-zero divisor in ${\calA}_+$. By Proposition (\ref{prop:discriminant}) the discriminant is $\alt^2(x)$, which by definition is a unit in ${\calA}$.
\end{proof}

\subsection{The Grothendieck-Deligne norm map}
\label{subsec:GD} Let $B\ra E $ be an extension of rings, where $E$ is free of rank $n$ as a $B$-module. The norm induced by the determinant map $\determinant : E \ra B$ corresponds to a $B$-algebra homomorphism $\sigma_E : \TS^n_BE \ra B$ (\cite{roby_lois_pol_mult},  \cite[6.3, p.180]{sga4_deligne_coh_supp_prop}, and \cite{iversen_lineardeterminants}). The algebra homomorphism $\sigma_E$ takes $y\otimes \cdots \otimes y$ to $\determinant (e\mapsto ye)$. 

\begin{lemma}
\label{lemma:GD}
Let $B\ra E$ be an extension of $A$-algebras, and let $f : R\ra E$ be an $A$-algebra homomorphism. Assume that the elements $f(x)=f(x_1), \ldots, f(x_n)$ form a $B$-module basis of $E$, and let $d_E \in B$ denote the discriminant of the the extension. 
\begin{enumerate}
\item[\upshape{(i)}] The extension of the canonical ideal $I_{R}$  by the induced map $\TS^n_AR \ra \TS^n_BE$ coincides with the canonical ideal $I_{E}$ and the ideal generated by $\alt^2(f(x_1), \ldots, f(x_n))$. 
\item[\upshape{(ii)}] The element $\alt^2(f(x))$ in $\TS^n_BE$ is mapped to the discriminant $d_E$ by the homomorphism $\sigma_E : \TS^n_BE \ra B$.
\item[\upshape{(iii)}]  We have the following commutative diagram of $A$-algebras and $A$-algebra homomorphisms
\eqbeg
\label{eq:GD}
\vcenter{\xymatrix{
\TS^n_AR \ar[r]^{\can} \ar[d]^{\can} & \TS^n_BE \ar[r]^{\sigma_E} & B \ar[d]^{\can} \\
{\calA} \ar[rr]^{\n_{E/B}} & & B_{d_E}, 
}}
\eqend
where $\n_{E/B}$ is the morphism of Theorem (\ref{thm:univ.et}).
\end{enumerate}
\end{lemma}

\begin{proof} As $f(x_1), \ldots, f(x_n)$ is a $B$-module basis of $E$ it follows from (\ref{subsec:linspan}) that any alternating tensor $\alt(z_1, \ldots, z_n)$ is in the $B$-module spanned by $\alt(f(x))$. In particular we see that $I_{R}\TS^n_BE=I_{E}=(\alt^2(f(x)))$, and the first statement of the lemma follows.

Let $z$ be an element of the free $B$-module $E$. We have (cf. \cite{iversen_lineardeterminants} p.9, Section 2.5) that the homomorphism $\sigma_E$ sends the polarized power sum $\p(z)$ to the trace of the multiplication map $e\mapsto ze$. Thus, the map $\sigma_E$ and $\n_{E/B}$ have the same action on polarized powers. Furthermore, as $\alt^2(f(x))=\det (\p(f(x_ix_j)))$ \eqref{eq:traceexp} and we have that $f(x_1), \ldots , f(x_n)$ form a $B$-module basis of $E$ it follows that  $\alt^2(f(x))$ is mapped to the discriminant $d_E$ of $B\ra E$. We have then proved the second statement of the lemma. We have also that the localization of the composite map $\sigma_E \circ \can : \TS^n_AR \ra B$ with respect to the element $\alt^2(x)$ is the morphism $\n_{E/B}$, proving the  third statement. 
\end{proof}

\subsection{Induced map of blow-up algebras}
Let $B\ra E$ be an extension as in  Lemma (\ref{lemma:GD}). From the $A$-algebra homomorphism $\TS^n_AR \ra B$ appearing as the top horizontal row in \eqref{eq:GD}, we obtain a graded homomorphism of blow-up algebras 
\eqbeg
\label{eq:degreezero}
\oplus _{m \geq 0}I^m_{R} \ra \oplus _{m \geq 0}I^m_{R}B=\oplus_{m \geq 0}(d_E)^m.
\eqend
We have furthermore by Lemma (\ref{lemma:GD}) that $\alt^2(x)$ in $\TS^n_AR$ is mapped to the discriminant $d_E$ by the map \eqref{eq:degreezero}. If we let $B_+$ denote the degree zero part of the localization of $\oplus_{m\geq0}I_{R}^mB$ with respect to the degree one element $d_E$, we obtain from \eqref{eq:degreezero} an induced $A$-algebra homomorphism
\eqbeg
\label{eq:n^+}
 \n^+_{E/B} : {\calA}_+ \ra B_+.
\eqend
Let $E_+:= E\otimes_B B_+$.

We will in the sequel use the notation  $\ann (d_E)\subseteq B$ for the kernel of the localization map $B \ra B_{d_E}$. With this notation we have
\eqbeg
\label{eq:annihilator}
B_+=B/\ann(d_E),
\eqend
as $I_{R}B=(d_E)$.

\begin{lemma}
\label{lemma:free} 
Let $A\ra R$ be a homomorphism of rings, and assume that $x_1, \ldots, x_n$ form an $A$-module basis for $R$. Then the induced $A$-algebra homomorphism  $\n_{R/A}^+ : {\calA}_+ \ra A_+$ \eqref{eq:n^+} is an isomorphism, and the induced map ${\calR}_+ \ra R_+$ is an isomorphism.
\end{lemma}

\begin{proof} We have that $A_+=A/\ann(d_R)$, where $d_R$ is the discriminant of $A\ra R$. Since $\n^+_{R/A} : {\calA}_+ \ra A_+$ is an $A$-algebra homomorphism it is necessarily surjective. Injectivity we prove in the following way. From Lemma (\ref{lemma:GD}) we have the following commutative diagram of $A$-algebras
\eqbeg
\label{eq:locdiagram}
\vcenter{\xymatrix{
{\calA}_+ \ar[d] \ar[r]^{\n^+_{R/A}} & A_+ \ar[d] \\
{\calA} \ar[r]^{\n_{R/A}} & A_{d_R}.
}}
\eqend
The vertical arrows are injective being the localizations in their respective discriminants,  and the bottom horizontal map $\n_{R/A}$ is an isomorphism by Proposition (\ref{prop:whenRisfree}). Injectivity of $\n_{R/A}^+$ now follows from the commutative diagram \eqref{eq:locdiagram}. 

As ${\calR}_+$ is both an ${\calA}_+$-algebra and an $R$-algebra we have an induced $A$-algebra homomorphism
\eqbeg
\label{eq:again}
R\otimes_A{\calA}_+ \ra {\calR}_+.
\eqend
The map  (\ref{eq:again}) sends the  ${\calA}_+$-module basis $x_1\otimes 1,\ldots, x_n\otimes 1$ to $\varphi_n(x_1),\ldots, \varphi_n(x_n)$.  By Proposition (\ref{prop:specialcase}) it follows that the map \eqref{eq:again} is an isomorphism of ${\calA}_+$-modules, hence an isomorphism of algebras. As $R_+$ is by definition $R\otimes_AA_+$, and we have $A_+={\calA}_+$ we obtain from \eqref{eq:again} the isomorphism ${\calR}_+ \ra R_+$.
\end{proof}

\begin{thm}
\label{thm:existence} Let $A\ra R$ be a homomorphism of algebras, and let $x_1, \ldots , x_n$ be elements of $R$. Let $B\ra E $ be an extension of $A$-algebras, and  assume that $f : R\ra E$ is an $A$-algebra homomorphism such that the elements $f(x_1), \ldots, f(x_n)$ form a $B$-module basis of $E$. Then the $A$-algebra homomorphism $\n^+_{E/B} : {\calA}_+ \ra B_+$ \eqref{eq:n^+} is the unique $A$-algebra homomorphism such that 
$$ {\calR}_+\otimes_{{\calA}_+}B_+ = E_+$$
as quotients of $R\otimes_{A}B_+$.
\end{thm}

\begin{proof} From the canonical map $\TS^n_AR \ra \TS^n_BE$ we obtain an induced $A$-algebra homomorphism of graded rings
$$  \oplus_{m\geq 0}I_{R}^m \ra \oplus_{m\geq 0}I_{R}^m\TS^n_BE.$$
By Lemma (\ref{lemma:GD}) (i) we have that the extension $I_R\TS^n_BE$ is the canonical ideal $I_E$. We let $({\calA}_E)_+$ denote the degree zero part of the localization of the graded ring $\oplus_{m\geq 0}I_{R}^m\TS^n_BE=\oplus_{m\geq0}I_E^m$ at the degree one element $\alt^2(f(x))$. And similarly we let $({\calR}_E)_+$ denote the localization of the graded ring we obtain from the natural map $\TS^n_BE \ra \TS^{n-1,1}_BE$. From the commutative diagram of $A$-algebras in \eqref{eq:commdiag} we obtain the following commutative diagram
\eqbeg
\label{eq:commdiag2}
\vcenter{\xymatrix{
{\calA}_+ \ar[d] \ar[r]& ({\calA}_E)_+ \ar[d] \\
{\calR}_+ \ar[r] & ({\calR}_E)_+
}}
\eqend 
As $({\calA}_E)_+$ is the localization of the blow-up algebra of the canonical ideal $I_E\subseteq \TS^n_BE$ we have by Proposition (\ref{prop:specialcase}) that the elements $\varphi_n(f(x_1)), \ldots, \varphi_n(f(x_n))$ form an $({\calA}_E)_+$-module basis for the algebra $({\calR}_E)_+$. The commutative diagram \eqref{eq:commdiag2} induces a canonical ${\calA}_+$-algebra homomorphism
\eqbeg
\label{eq:ind}
{\calR}_+ \otimes_{{\calA}_+} ({\calA}_E)_+ \ra ({\calR}_E)_+.
\eqend
The ${\calA}_+$-algebra homomorphism \eqref{eq:ind} identifies the $({\calA}_E)_+$-module basis $\varphi_n(x_1)\otimes 1, \ldots ,\varphi_n(x_n)\otimes 1$ with the basis $\varphi_n(f(x_1)), \ldots, \varphi_n(f(x_n))$. And consequently \eqref{eq:ind} is an isomorphism of ${\calA}_+$-algebras. Furthermore, by Lemma (\ref{lemma:free}) we have a natural identification $({\calA}_E)_+ =B_+$ and $({\calR}_E)_+=E_+$. Hence we have shown the existence of a morphism ${\calA}_+ \ra B_+$ with the desired property.

For uniqueness we note that a morphism $\eta : {\calA}_+ \ra B_+$ such that ${\calR}_+\otimes_{{\calA}_+}B_+=E_+$ would have to map the discriminant $d^+$ of ${\calA}_+ \ra {\calR}_+$ to the discriminant $d_E$ of $B \ra E$. When we localize ${\calA}_+$ in the discriminant $d^+$ we obtain ${\calA}$, and consequently a commutative diagram of $A$-algebras
\eqbeg
\label{eq:unique}
\vcenter{\xymatrix{
{\calA}_+ \ar[r]^{\eta} \ar[d]^{\can} & B_+ \ar[d]^{\can} \\
{\calA} \ar[r] &  B_{d_E}. 
}}
\eqend

As $B_{d_E} \ra E_{d_E}$ is \'etale it follows from Theorem (\ref{thm:univ.et}) that the bottom horizontal row of the above \eqref{eq:unique} diagram is $\n_{E_{d}/B_{d}}$. Consequently the morphisms $\eta $ and $\n_{E/B}^+$ coincide after localization with respect to the discriminant $d^+$. As the vertical arrows in the diagram above are injective, both discriminants being non-zero divisors, it follows that $\eta =\n_{E/B}^+$.
\end{proof}

\subsection{Generically \'etale families}
Let $B\ra E$ be a finite and locally free homomorphism of rings. If the local generators of the discriminant ideal $D_{E/B}$ are non-zero divisors then we say that the extension $B\ra E$ is {\em generically \'etale}. Equivalently, the family is generically \'etale if the closed subscheme of $\Spec (B)$, corresponding to the ideal $D_{E/B}$, is an effective Cartier divisor. Another, perhaps more geometric, characterization is the following. The finite and locally free extension $B\ra E$ is generically \'etale if and only if the open subset $U\subseteq \Spec (B)$ of points where the fibers are \'etale is schematically dense in $\Spec (B)$.

\subsection{The universal generically \'etale family, local description}\label{subsec:genetale}
In analogy with the end of Section 3 we give here a local treatment for the universal properties of the space of generically \'etale families. Let ${\bf x}$ denote the set of all $n$-tuples of elements $x_1, \ldots, x_n$ in $R$. We let
$$ \calU_{R}^+ =\bigcup_{x\in {\bf x}} \Spec (\calA_+(\alt^2(x))) \subseteq \Proj (\oplus_{m\geq 0} I_{R}^m).$$
Using identities similar to those in (\ref{subsec:nu^2(x)cover}) one sees that the open subset $\calU_{R}^+$ equals the scheme $\Proj(\oplus_{m\geq 0} I_{R}^m)$. We note furthermore, that the family $\calZ_{R}^+=\cup_{x\in {\bf x}} \Spec (\calR_+(\alt^2(x))) \ra \calU_{R}^+$ is generically \'etale, and is a closed subscheme of $\Spec (R)\times_{\Spec (A)}\calU_{R}^+$. 

\begin{cor}
\label{cor:gen.et} 
Let $B$ be an $A$-algebra, and let $E$ be an algebra quotient of $R\otimes_AB$. Assume that the homomorphism of rings $B\ra E$ is a generically \'etale extension of rank $n$. Then there exist a unique morphism of schemes
$$ \n ^+ : \Spec (B) \ra \calU_{R}^+,$$
such that the pull-back of $\calZ_{R}^+$ along $\n ^{+}$ is $\Spec (E)$. 
\end{cor}

\begin{proof}The existence and uniqueness of the homomorphism $\n^+$ follows by Theorem (\ref{thm:existence}) and arguments similar as in the proof of Corollary (\ref{cor:etalaffine}).
\end{proof}


\section{The space of \'etale families}

In this section we show how to construct a space parameterizing \'etale families in a fixed separated algebraic space $X\ra S$ (\cite[6.6]{laumonstacks}). We thereby obtain an abstract and global version of Corollary (\ref{cor:etalaffine}), but not the explicit statement about the basis as in Theorem (\ref{thm:univ.et}). 

\subsection{Disjoint sections}
Let $f : X \ra S$ be a morphism of algebraic spaces. A {\em section} of $f$ is a morphism $s : S \ra X$ such that $f\circ s =\operatorname{id}_S$. Two sections $s, s' : S \ra X$ are called {\em disjoint} if $(s,s') : S \ra X\times_SX$ does not intersect the diagonal. 

\begin{rem} For a separated morphism $X\ra S$ we have that a section $s : S \ra X$ is equivalent with a closed subspace $Z\subseteq X$, with $Z$ isomorphic to $S$. Note that two disjoint sections $s, s' : S \ra X$ determines a closed subspace $Z$ isomorphic to the disjoint union of two copies of $S$.
\end{rem}

 For any $S$-space $T$ we let $U^n_X(T)$ denote the set of {\em unordered} $n$-tuples of sections $s_1, \ldots ,s_n$ of the second projection $X\times_ST\ra T$, where the sections are pairwise disjoint. Clearly $U^n_X$ describes a functor; the functor of $n$ {\em unordered disjoint sections} of $X\ra S$.


\subsection{Diagonals and their complement}
Let $X \ra S$ be separated, and let $X^n_S:=X\times_S \cdots \times_SX$ denote the $n$-fold product. We let $\cdiag$ denote the open complement of the diagonals in $X^n_S$.

\subsection{Action of the symmetric group}
The permutation action $\Sym_n$ on $X_S^n$ induces an action $\rho : \cdiag \times \Sym_n \ra \cdiag$. It is clear that the map $(\pi_1, \rho) : \cdiag \times \Sym_n \ra \cdiag \times_S \cdiag$ is a monomorphism of algebraic  spaces, where $\pi_1 $ denotes the projection on the first factor. We say that the action $\rho $ is {\em free}, and we note that $(\pi_1, \rho)$ describes an \'etale equivalence relation on $\cdiag$.

\begin{lemma}
\label{lemma:presheaf} Let $X \ra S$ be a separated algebraic space. Then $U^n_X$ is the presheaf quotient of the equivalence relation $\xymatrix@M=1pt{ {\cdiag\times\Sym_n}
\ar@<.5ex>[r]^-{\pi_1} \ar@<-.5ex>[r]_-{\rho} & {\cdiag}}$.
\end{lemma}
\begin{proof} Clearly $X_S^n$ parameterizes $n$ ordered sections as we have the canonical identification
$$ \hom_S(T, X_S^n)=\prod_{i=1}^n \hom_T(T,X\times_ST),$$
for any $S$-space $T$. Similarly we get that $\cdiag$ parameterizes $n$ ordered disjoint sections of $X \ra S$. The action of $\Sym_n$ on $\cdiag$ corresponds to the action on $\hom_S(T,X_S^n)$ given by permuting the sections, from where it follows that the presheaf quotient is $U^n_X$.
\end{proof}

\subsection{The space of \'etale families}\label{subsec:sheafU}
We let ${\calU}^n_X$ denote the sheafification of the presheaf $U^n_X$ in the \'etale topology. As a consequence of the above lemma we have that the sheaf ${\calU}^n_X$ is the quotient sheaf of the equivalence relation $\xymatrix@M=1pt{ {\cdiag\times\Sym_n}
\ar@<.5ex>[r]^-{\pi_1} \ar@<-.5ex>[r]_-{\rho} & {\cdiag}}$.

\begin{prop}
\label{prop:sheafU}
Let $X\ra S$ be a separated algebraic space. The sheaf ${\calU}^n_X$ is an algebraic space and represents the functor of closed subspaces of $X$ that are \'etale and of rank $n$ over the base. 
\end{prop}

\begin{proof} As the algebraic space $V_X\ra S$ is separated, the map $(\pi_1, \rho) : \cdiag \times \Sym_n \ra \cdiag \times_S \cdiag$ is a closed immersion and hence quasi-affine. We then have that the quotient sheaf  ${\calU}^n_X$  of the \'etale equivalence relation $\xymatrix@M=1pt{ {\cdiag\times\Sym_n}
\ar@<.5ex>[r]^-{\pi_1} \ar@<-.5ex>[r]_-{\rho} & {\cdiag}}$ is an algebraic space ~\cite[II, Prop. 3.14]{knutson_alg_spaces}.

What remains to show is that sections of ${\calU}^n_X$ corresponds to closed subspaces that are \'etale of rank $n$ over the base. By Lemma (\ref{lemma:presheaf}) a section over an $S$-space $T$ of the presheaf $U^n_X$ is $n$ unordered disjoint sections of $X\times_ST \ra T$. Equivalently, as $ X \ra S$ is separated, a section over $T$ is a closed subspace $Z\subseteq X\times_ST$ such that $Z$ is isomorphic to a disjoint union of $n$-copies of $T$. Consequently, a  section over $T$ of the sheaf ${\calU}^n_X$ is  given by an \'etale covering $T' \ra T$ and a closed subspace $Z' \subseteq X\times_ST'$, isomorphic to a disjoint union of $n$ copies of $T'$. The pull-back of the family $Z' \ra T'$ along the two different projection maps $T'\times_TT'\ra T'$ coincide. Therefore we have that the section over $T$ of ${\calU}_X^n$ is given by a unique closed subspace $Z\subseteq X\times_S T$ such that $Z\ra T$ is \'etale and finite of rank $n$.

Further, as every finite \'etale morphism $Z\ra T$ trivializes after an \'etale base
change, the set of sections over $T$ is the set of all  closed subspaces $Z\subseteq X\times_ST$ that are finite, \'etale and of rank $n$.
\end{proof}

\begin{rem} Let $\Sym_{n-1}$ act by permuting the first $(n-1)$ factors of $X^n_S$, and consider the induced action on $\cdiag$. The induced map $\cdiag /\Sym_{n-1} \ra \calU_X^n$ will be \'etale of rank $n$, and one can check that this is in fact the universal family of \'etale, rank $n$ closed subspaces of $X$. In Proposition (\ref{prop:implicit}) we will prove this fact using a different approach.
\end{rem}

\begin{rem}
\label{rem:support} When $X=\Spec (R)$ is affine over base $S=\Spec(A)$ then we have that ${\calU}^n_X$ is the open complement of the diagonals in $\Spec (\TS^n_AR)$. One can furthermore check that the support of the closed subscheme defined by the canonical ideal $I_{R}\subseteq \TS^n_AR$ is the diagonals. Hence Corollary (\ref{cor:etalaffine}) can be obtained as a consequence of Proposition (\ref{prop:sheafU}). However, the more explicit description given by Theorem (\ref{thm:univ.et}) is what  will be important for us.
\end{rem}

\section{The space of generically \'etale families}

\subsection{Finite group quotients}
When a finite group $G$ acts on a separated algebraic space $X$, the geometric quotient $X/G$ exists as an algebraic space; this is an unpublished result of Deligne ~\cite[p.183]{knutson_alg_spaces}. When the base is locally Noetherian and $X\ra S$ is locally of finite type, proofs of this  existence result are given in \cite{keel_mori_quotients} and in \cite{kollar_quotients}. It is furthermore possible to extend the proof given by Koll\'ar in  \cite{kollar_quotients} to the general setting with any separated algebraic space $X\ra S$, for details we refer to \cite{rydh_finite_quotients}.

 We will be interested in the particular case with the symmetric group $\Sym_n$ of $n$ letters acting by permuting the factors of the $n$-fold product $X^n_S=X\times_S \cdots \times_S X$. For a separated algebraic space $X\ra S$ we  will denote the quotient space with $\symquot ^n_SX$.

\subsection{Symmetric spaces}
\label{subsec:quot-flatbasechange}A well-known fact is that geometric quotients commute with flat base change \cite[p.9]{GIT} or, for finite groups \cite[Exp.V, Prop.1.9]{sga1}. Thus if $T \ra S$ is flat, and $X_T$ denotes $X\times_ST$, then we have the cartesian diagram
\eqbeg
\label{eq:cartsym}
\vcenter{\xymatrix{
\symquot^n_{T}(X_{T}) \ar[d] \ar[r] & \symquot^n_SX \ar[d] \\
T \ar[r] & S.}}
\eqend
This fact will be useful when we want to reduce to the case with an affine base scheme $S$.

\subsection{Fixed-point reflecting morphisms}\label{fpr}
\label{subsec:fprconstruction} We recall some properties of geometric quotients (see \cite{knutson_alg_spaces}, \cite{rydh_finite_quotients}). Let $G$ be a finite group acting on an algebraic  space $Y$. The stabilizer group $G_y$ of a point $y\in Y$ is defined as the inverse image of the point $(y,y)$ by the map $G\times_SY \ra Y\times_SY$. A $G$-equivariant  morphism $f : Y \ra X$ is {\em fixed-point reflecting} in a point $y\in Y$ if the stabilizer group $G_y$ equals the stabilizer group $G_{f(y)}$. The fixed-point reflecting set, with respect to a given \'etale separated map $Y\ra X$, we denote by $Y|_{\fpr}$. The set $Y|_{\fpr}\subseteq Y$ is an open $G$-invariant subset, and we denote the geometric quotient $Y|_{\fpr}/G $ by $Y/G|_{\fpr}$. We have furthermore a cartesian diagram
\eqbeg
\label{eq:fprcartesian}
\vcenter{\xymatrix{
Y| _{\fpr} \ar[d] \ar[r] & X \ar[d] \\
Y/G|_{\fpr} \ar[r] & X/G,}}
\eqend
where the horizontal maps are \'etale.

\begin{prop}\label{prop:closedimmersion}
Let $X \ra S$ be a separated algebraic space, and let $Y \ra X$ be \'etale, with $Y$ affine. Let $T \ra S$ be an $S$-space, and let $Z\subseteq Y\times_ST$ be a closed subspace, which is finite and locally free of rank $n$ over $T$. Let $\sigma_Z : T \ra \symquot^n_SY$ denote the morphism corresponding to the determinant map (\ref{subsec:GD}). Then the morphism $\sigma_Z : T \ra \symquot^n_SY$ factors through the fixed-point reflecting set $\symquot^n_SY| _{\fpr}$ if and only if the composite morphism $Z \subseteq Y\times_ST \ra X\times_ST$ is a closed immersion.
\end{prop}

\begin{proof} Since the morphism $Z\ra T$ is finite we have that the morphism $i : Z \ra X\times_ST$ is a closed immersion if and only if the induced morphism over points is a closed immersion. Thus we can assume $T=\Spec (K)$, where $K$ is an algebraically closed field. The support of the family $Z$ is then a finite set of points $z_1, \ldots , z_p$ in $Y \times_S \Spec (K)$. Let $m_i$ denote the length of the local ring at $z_i$, $i=1, \ldots , p$. The determinant map $T \ra \symquot^n_TZ$ takes the family $Z \ra T$ to the cycle $m_1\cdot z_1 + \cdots + m_p\cdot z_p$ ~\cite[Prop.4.7]{iversen_lineardeterminants}. The determinant map composed with the induced map $\symquot^n_TZ \ra \symquot^n_SY$ is the morphism $\sigma_Z$. Let $z\in Y^n_S$ be the $n$-tuple of points  where the $m_1$ first coordinates are $z_1$, the next $m_2$ coordinates are $z_2$ etc. The set of points in $Y^n_S$ that we obtain by shuffling the order of the $n$-tuple $z$ is the inverse image $q_n^{-1}(\sigma_Z(T))$ by the quotient map $q_n : Y^n_S \ra \symquot^n_SY$.

Since $Z\subseteq Y\times_ST$ is a closed immersion and $Y\ra X$ an \'etale map, the composition $i : Z \ra X\times_ST$ is a closed immersion if and only if it is injective on points (\cite[Prop. 17.2.6, Cor. 18.12.6]{egaIV}).  The map $i$ being  injective is equivalent with the set $q_n^{-1}(\sigma_Z(T))\subset Y^n\times_S T$ being a subset of the fixed-point reflecting set $(Y^n_S)| _{\fpr}$ with respect to the map $Y^n_S \ra X^n_S$.
\end{proof}

\subsection{Covers of symmetric quotients with affine base} 
Let the base space $S$ be affine, and let $X \ra S$ be a separated algebraic space. We choose an \'etale covering $\coprod X_{\gamma} \ra X$, with $X_{\gamma}$ affine for each $\gamma$, such that the induced map 
\eqbeg
\label{eq:alphacover}
 \coprod_{\gamma} (X_{\gamma})^n_S \ra X^n_S
\eqend
is surjective. In other words, a covering such that any $n$-tuple of points in $X$ lies in the image of some $X_{\gamma} \ra X$. Then $\coprod_\gamma\symquot^n_S(X_\gamma)|_{\fpr} \ra \symquot^n_S X$ is an \'etale cover. For two indices $\gamma$ and $\gamma'$ we define
$$(X_{\gamma}\times_X X_{\gamma'})^n| _{\fpr}:= (X_{\gamma})^n| _{\fpr} \times_{X^n_S} (X_{\gamma'})^n| _{\fpr},$$
where we have suppressed the base $S$ in the notation, and where we have used the notation introduced in (\ref{subsec:fprconstruction}). We have the following commutative diagram
\eqbeg
\label{eq:opensym}
\vcenter{\xymatrix{
 {\coprod_{\gamma, \gamma'} (X_{\gamma}\times_X X_{\gamma'})^n| _{\fpr}} \ar[d]^q \ar@<.5ex>[r]^-{p_1} \ar@<-.5ex>[r]_-{p_2} & \coprod_{\gamma}(X_{\gamma})^n| _{\fpr}  \ar[d]^q \ar[r]^-p & X^n_S \ar[d]^q \\
{\coprod _{\gamma, \gamma'} \symquot^n_S(X_{\gamma}\times_X X_{\gamma'})| _{\fpr}}\ar@<.5ex>[r]^-{\pi_1} \ar@<-.5ex>[r]_-{\pi_2} & \coprod_{\gamma} \symquot^n_S(X_{\gamma})| _{\fpr} \ar[r]^-{\pi} & \symquot^n_SX,  }}
\eqend
where the $q$'s are the quotient maps, $p_1, p_2$ and $p$ are the natural maps, and where $\pi_1, \pi_2$ and $\pi$ are the induced ones.  By (\ref{fpr}) we have that the three squares: $qp=\pi q$ and $qp_i=\pi_iq $ for $i=1, 2$; are cartesian. Furthermore we have that $\pi_1$ and $\pi_2$ form an \'etale equivalence relation with quotient $\pi$.

\begin{lemma} \label{lemma:normideal-basechange}
Let $A\ra A'$ be a homomorphism of rings, $R$ be an $A$-algebra and $R'=R\otimes_A A'$. Then the extension of the canonical ideal $I_R$ (\ref{subsec:normideal}) by the homomorphism $\TS^n_A(R)\ra\TS^n_{A'}(R')$ equals the canonical ideal $I_{R'}$.
\end{lemma}
\begin{proof}
Recall the alternator map $\alt_R : \Ten^n_AR \ra  \Ten^n_AR$ \eqref{eq:normmap} and that the canonical ideal $I_R$ is generated by $\alt(x)\alt(y)$ for $x,y\in \Ten^n_AR$. As $\alt_{R'}=\alt_R\otimes 1 : \Ten^n_{A'}R' \ra \Ten^n_{A'}R'$ is $A'$-linear the lemma follows.
\end{proof}

\begin{lemma} \label{lemma:normideal-etalefpr}
Let $X'\ra X$ be an \'etale morphism of $S$-schemes where $S=\Spec(A)$, $X'=\Spec(R')$ and $X=\Spec(R)$ are affine schemes. Let $\calI$ and $\calI'$ be the ideal sheaves corresponding to the canonical ideals of $\TS^n_A(R)$ and $\TS^n_A(R')$ respectively. Then the pull-back of $\calI$ by the morphism
$$\symquot^n_S(X')|_\fpr \ra \symquot^n_S(X)$$
equals the restriction of $\calI'$.
\end{lemma}

\begin{proof}
The canonical ideal $I_R$ of $\TS^n_A R$ is the image of the $\TS^n_AR$-linear map
\eqbeg
\label{eq:I_R}
\alt\times \alt : \Ten^n_AR \otimes_{\TS^n_AR}\Ten^n_AR \ra \TS^n_AR.
\eqend
Let $f_\gamma\in \TS^n_A(R')$ be symmetric tensors such that the $D(f_\gamma)$ covers $\symquot^n_S(X')|_\fpr$. 
From \eqref{eq:fprcartesian} we have that the diagram
$$\xymatrix{
(\Ten^n_AR')_{f_{\gamma}} & \ar[l] \Ten^n_AR \\
(\TS^n_AR')_{f_{\gamma}} \ar[u] &\ar[l]  \TS^n_AR \ar[u],}
$$
is co-cartesian for every $\gamma$. Therefore, by applying the change of basis $\TS^n_AR \ra (\TS^n_AR')_{f_\gamma}$ to the morphism \eqref{eq:I_R} we obtain
$$\alt\times \alt : (\Ten^n_AR')_{f_\gamma} \otimes_{(\TS^n_AR')_{f_\gamma}}(\Ten^n_AR')_{f_\gamma} \ra (\TS^n_AR')_{f_\gamma}.$$
Thus the extension of the canonical ideal $I_R$ in $(\TS^n_AR')_{f_\gamma}$ equals the localization of the canonical ideal $I_{R'}$ in $f_\gamma$.
\end{proof}

\begin{prop} \label{prop:normideal-glues} Let $X\ra S$ be a separated algebraic space. The canonical ideals $I_{R}\subseteq \TS^n_AR$ as defined in Section (\ref{subsec:normideal}) for $A$-algebras $R$, glue to an ideal sheaf ${\calI}_{X}$ on $\symquot^n_S(X)$. If $Y\ra X$ is an \'etale morphism then the pull-back of ${\calI}_{X}$ along $\symquot^n_S(Y)|_\fpr\ra\symquot^n_S(X)$ is the restriction of the canonical ideal sheaf ${\calI}_Y$ on $\symquot^n_S(Y)$.
\end{prop}

\begin{proof} Let $\coprod_\beta S_\beta \ra S$ be an \'etale cover with $S_\beta$ affine and let $X_\beta=X\times_S S_\beta$. For any $\beta$ choose an \'etale cover $\coprod_{\gamma} X_{\beta,\gamma}\ra X_\beta$ as in \eqref{eq:alphacover}. Then the canonical ideal sheaves on $\symquot_{S_\beta}(X_{\beta,\gamma})$ glues to an ideal sheaf $\calI_{X_\beta}$ on $\symquot_{S_\beta}(X_\beta)$ using Lemma~(\ref{lemma:normideal-etalefpr}) and the \'etale equivalence relation of \eqref{eq:opensym}. It is further clear that $\calI_{X_\beta}$ is independent on the choice of \'etale cover of $X_\beta$.

Using that symmetric products commute with flat base change \eqref{eq:cartsym} and Lemma~(\ref{lemma:normideal-basechange}) it follows that the sheaves $\calI_{X_\beta}$ glue to an ideal sheaf ${\calI}_{X}$ on $\symquot^n_S(X)$ which is independent on the choice of covering of $S$. The last statement is obvious from the construction.
\end{proof}

\subsection{The addition of points map}
Let $q_n : X_S^n \ra \symquot^n_SX$ be the quotient map. We let $\Sym_{n-1}$ act on the first $(n-1)$ copies of $X_S^n$ and denote the geometric quotient
$X_S^n/\Sym_{n-1}$ with $\symquot^{n-1,1}_SX$. Since the geometric quotient equals the categorical quotient in the category of separated algebraic spaces (\cite[Cor. 2.15]{kollar_quotients}, \cite{rydh_finite_quotients}), and as $q_n$ is $\Sym_{n-1}$-invariant it factors through $X_S^n \ra \symquot^{n-1,1}_SX$. We have a canonical morphism
\eqbeg
\label{eq:psi-fam}
 \psi_X : \symquot^{n-1,1}_SX \ra \symquot^n_SX.
\eqend

\begin{lemma}\label{lem:psibehaveswell}
Let $Y\ra X$ be an \'etale map. Then the pull-back of the family $\psi_X$ along the morphism $\symquot^n_SY|_{\fpr} \ra \symquot^n_SX$ is canonically identified with the restriction of the family $\psi_Y$ to the open subset $\symquot^n_SY|_{\fpr} \subseteq \symquot^n_SY$.
\end{lemma}

\begin{proof} Let $U=(Y^n_S)|_{\fpr}$ denote the $\Sym_n$-invariant fixed-point reflecting set of $Y^n_S\ra X^n_S$. We then have that the induced map $U/\Sym_{n-1} \ra U/\Sym_n$ is the restriction of $\psi_Y : Y^n/\Sym_{n-1}\ra Y^n/\Sym_n$ to the open subset $ U/\Sym_n=\symquot^n_S(Y)|_{\fpr}\subseteq \symquot^n_SY$. We need to see that the pull-back of the family $\psi_X : X^n/\Sym_{n-1} \ra X^n/\Sym_n$ along the map $U/\Sym_n \ra X^n/\Sym_n$ equals $U/\Sym_{n-1} \ra U/\Sym_n$. First note that $U$ is the base change of $X^n$ by the \'etale morphism $U/\Sym_n \ra X^n/\Sym_n$ as the diagram \eqref{eq:fprcartesian} is cartesian. As taking the quotient with $\Sym_{n-1}$ commutes with flat base change (\ref{subsec:quot-flatbasechange}) the result follows.
\end{proof}

\begin{lemma}
\label{lem:lifting} Let $Z \ra S$ be a finite map of algebraic spaces, and let $W\ra Z$ be an \'etale cover. Then there exists an \'etale cover $S'\ra S$ such that $W\times_SS' \ra Z\times_SS'$ has a section. In particular when $Z\subseteq X$ is a closed immersion of algebraic spaces, with $Z$ finite over the base. Then, for a separated and \'etale cover $Y \ra X$ there exists an \'etale cover $S'\ra S$ such that the closed immersion $Z\times_SS' \ra X\times_SS'$ lifts to a closed immersion $Z\times_SS' \ra Y\times_SS'$.
\end{lemma}

\begin{proof} Let $x \in S$ be a point, and let $A^{h}_x$ denote the strictly local ring at $x$. Let $E_x$ be coordinate ring of the affine scheme $Z\times_S\Spec(A^{h}_x)$. Then $E_x$ is a product of local Henselian rings with separably closed residue fields. For every closed point $z_i$ of $\Spec(E_x)$ choose a point $w_i$ of $W\times_S\Spec(A^{h}_x)$ above $z_i$. As $W\ra Z$ is \'etale there is a section $Z\times_S\Spec(A^{h}_x))\ra W\times_S\Spec(A^{h}_x)$ mapping $z_i$ to $w_i$. By a standard limit argument, this section extends to a section $Z\times_S U \ra W\times_S U$ where $U$ is an \'etale neighborhood around the point $x \in S$, proving the first claim.
\end{proof}

\begin{prop}
\label{prop:implicit}
Let $X\ra S$ be a separated algebraic space. Let $\Delta\subseteq \symquot_S^nX$ denote the closed subspace defined by the canonical ideal sheaf ${\calI}_{X}$. We have a canonical identification  ${\calU}_X^n =\symquot^n_SX \setminus \Delta$, where ${\calU}_X^n$ is the algebraic space of (\ref{subsec:sheafU}). Moreover, the family $\psi^{-1}_X({\calU}^n_X) \ra {\calU}^n_X$ is the universal family of closed subspaces in $X\ra S$, that are \'etale of rank $n$ over the base.
\end{prop}

\begin{proof} Denote the open subspace $U=\symquot^n_SX\setminus\Delta$, and let $Z=\psi_X^{-1}(U)$. By Proposition (\ref{prop:sheafU}) it suffices to show that the pair $(U,Z)$ represents the functor parameterizing closed subspaces of $X$ that are \'etale of rank $n$ over the base.

We can by (\ref{subsec:quot-flatbasechange}) assume that the base $S=\Spec (A)$ is affine, and we let $\coprod_{\gamma}X_{\gamma} \ra X$ be an \'etale covering as given in the construction \eqref{eq:alphacover}.  Let $U_{\gamma}$ denote the inverse image of $U$ along $\symquot^n_{S}(X_{\gamma})|_{\fpr} \ra \symquot^n_SX$. Let $Z_{\gamma}$ be the pull-back of the family $Z\ra U$ to $U_{\gamma}$ and define $U_{\gamma, \gamma'}$ and $Z_{\gamma, \gamma'}$ in a similar way. We then have a cartesian diagram of \'etale equivalence relations
\eqbeg\label{D:implicit}
\vcenter{\xymatrix{
{\coprod_{\gamma, \gamma'} Z_{\gamma, \gamma'} \ar[d]} \ar@<.5ex>[r]^-{p_1} \ar@<-.5ex>[r]_-{p_2} & \coprod_{\gamma} Z_{\gamma} \ar[d] \ar[r]^-p & Z \ar[d] \\
{\coprod _{\gamma, \gamma'} U_{\gamma, \gamma'}}\ar@<.5ex>[r]^-{\pi_1} \ar@<-.5ex>[r]_-{\pi_2} & \coprod_{\gamma} U_{\gamma} \ar[r]^-{\pi} & U.  }}
\eqend
By Proposition (\ref{prop:normideal-glues}) we have that $U_{\gamma}=\symquot^n_S(X_{\gamma})|_{\fpr}\cap D({\calI}_{X_{\gamma}})$, where $D({\calI}_{X_{\gamma}})$ is the open subset of $\symquot^n_S(X_{\gamma})$ defined by the non-vanishing of the ideal sheaf ${\calI}_{X_{\gamma}}$. By Lemma (\ref{lem:psibehaveswell}) the family $Z_{\gamma} \ra U_{\gamma}$ is simply the restriction of the family $\psi_{X_{\gamma}} $ to $U_{\gamma}$. Now it follows from Corollary (\ref{cor:etalaffine}) and Proposition (\ref{prop:closedimmersion}) that  $(U_{\gamma},Z_{\gamma})$ parameterizes closed subspaces $W\subseteq X_{\gamma}$ that are \'etale of rank $n$ over the base, and which are also closed  subspaces of $X$. In particular, we have that $Z\to U$ is \'etale of rank $n$. The universal properties of $(U,Z)$ then follows from Lemma (\ref{lem:lifting}) and the above diagram.
\end{proof}

\begin{rem}\label{rem:diagonals} It follows from the proposition that the support $|\Delta |$ of the space defined by the canonical sheaf of ideals  ${\calI}_{X}$,  equals the diagonals. This can also be verified directly.

\end{rem}
\subsection{The space of generically \'etale families}
\label{subsec:construction}
We let ${\calG}^n_X$ denote the blow-up of $\symquot^n_SX$ along the closed subspace $\Delta\subseteq \symquot^n_SX$ defined by canonical ideal ${\calI}_X$. We furthermore let ${\calZ}$ denote the blow-up of $\symquot^{n-1,1}_S(X)$ along $\psi^{-1}_X{\Delta}$, where $\psi_X $ is the canonical morphism \eqref{eq:psi-fam}. The morphism $\psi_X$ then induces a morphism
\eqbeg
\label{eq:family} 
{\calZ}_X \ra {\calG}^n_X.
\eqend

\begin{rem} The property of being generically \'etale is not stable under base change, as one easily realizes by taking the fiber of a point in the discriminant locus of the family. However the property of being generically \'etale is stable under flat base change, and in particular under \'etale base change.
\end{rem}

\begin{thm}
\label{thm:main} 
Let $X \ra S$ be a separated algebraic space, and let ${\calZ}_X\ra {\calG}^n_X$ be as in \eqref{eq:family}. Then the family ${\calZ}_X \ra {\calG}^n_X$ is generically \'etale of rank $n$, and has the following universal property. For any $S$-space $T$, and any closed subspace $Z\subseteq X\times_ST$ such that the projection $Z \ra T$ is generically \'etale of rank $n$, there exists a unique morphism $f : T\ra {\calG}^n_X$ such that the pull-back $f^{*}{\calZ}_X=Z$, as subspaces of $X\times_ST$.
\end{thm}  

\begin{proof} Proceeding as in the proof of Proposition (\ref{prop:implicit}), replacing $Z\ra U$ with $\calZ_X\ra\calG^n_X$, we obtain a cartesian diagram similar to \eqref{D:implicit}. In this diagram, the vertical arrows are the blow-ups of the canonical morphisms $\psi_{X_\gamma}$ and $\psi_{X_{\gamma,\gamma'}}$, in the corresponding canonical ideals, restricted to the fixed-point reflecting loci. This is because blowing up commutes with flat base change. Arguing as in the proof of Proposition (\ref{prop:implicit}), it then follows that ${\calZ}_X \ra {\calG}^n_X$ is generically \'etale and has the ascribed universal property by (\ref{subsec:genetale}), Corollary (\ref{cor:gen.et}), Proposition (\ref{prop:closedimmersion}) and Lemma (\ref{lem:lifting}).
\end{proof}

\subsection{Schematic closure}
Let $f : Y\ra X$ be a {\em quasi-compact} immersion of algebraic spaces. Then  $f_*{\calO}_Y$ is quasi-coherent, and in particular the kernel of ${\calO}_X \ra f_*{\calO}_Y$ is a quasi-coherent sheaf of ideals in ${\calO}_X$ ~\cite[II. Prop. 4.6]{knutson_alg_spaces}. The closed subspace in $X$ determined by this ideal sheaf is the {\em schematic closure}  of $f : Y\ra X$. 

\begin{cor} Let ${\calH}^n_X$ denote the Hilbert functor of flat, finite rank $n$ families of closed subspaces in $X\ra S$. Then  ${\calU}^n_X$ is an open subspace of ${\calH}^n_X$, and ${\calG}^n_X$ equals the schematic closure of ${\calU}^n_X$ in ${\calH}^n_X$.
\end{cor}

\begin{proof} As $X\to S$ is separated the Hilbert functor is representable by an algebraic space (see e.g. \cite{rydh_hilbert}, \cite{ekedahl&skjelnes}). It is clear that ${\calU}^n_X$ is open in ${\calH}^n_X$ being the complement of the discriminant of the universal family. As the discriminant is a locally principal subspace we have that the open immersion ${\calU}^n_X \subseteq {\calH}^n_X$ is quasi-compact. Let $G\subseteq {\calH}^n_X$ denote the schematic closure  of ${\calU}^n_X$. It is clear that the restriction of the universal family $\xi\to {\calH}^n_X$ to $G$ satisfies the universal property for generically \'etale families. Consequently, by the theorem we have that ${\calG}^n_X=G$. 
\end{proof}

\begin{rem} Let $R$ be an $A$-algebra. There is a canonical $A$-algebra homomorphism $\tau : \Gamma^n_A(R)\ra \TS^n_AR$ from the divided powers algebra to the ring of symmetric tensors. The map $\tau$ is in general  neither surjective nor injective \cite{lundkvist}. 

The situation can be globalized (see \cite{rydh_reprofzerocycles}) for a separated algebraic space $X \ra S$, giving a map
$ t : \symquot^n_SX \ra \Gamma^n_{X/S}$.
In (\cite{ekedahl&skjelnes}) they defined a closed subspace $\Delta'\subseteq \Gamma^n_{X/S}$ whose blow-up yields the good component of the Hilbert functor ${\calH}^n_X$. The closed subspace $\Delta \subseteq \symquot^n_{S}X$ that we consider in this article is the inverse image $t^{-1}(\Delta')$. Even though the map $ t : \symquot^n_SX \ra \Gamma^n_{X/S}$ is not an isomorphism we have that the two corresponding blow-ups of $\Delta'$ and $t^{-1}(\Delta')$ are isomorphic; indeed both blow-ups are identified with the schematic closure of ${\calU}^n_X$ inside the Hilbert space ${\calH}^n_X$.
\end{rem}

\begin{rem} The universal family ${\calZ}_X$ we obtained by blowing up the closed subspace $\psi^{-1}_X(\Delta) \subseteq \symquot^{n-1,1}_SX$, or equivalently by taking the strict transform of $\psi_X$ along ${\calG}^n_X \ra \symquot^n_SX$.  There is however a natural way to make a finite family flat. Let $\mathrm{B}$ denote the blow up of $\symquot^n_SX$ along the $n$-th Fitting ideal of the family $\psi_X : \symquot^{n-1,1}_SX \ra \symquot^n_SX$. Then the strict transform $E\ra \mathrm{B}$ of the family $\psi_X$ is flat  (\cite[5.4]{raynaud-gruson}). The Fitting ideal is in general different from the canonical ideal.
\end{rem}



\bibliographystyle{dary}

\bibliography{genetalefams}

\end{document}